\documentclass[11pt, twoside, fleqno]{article} %%, draft

\usepackage{stmaryrd}
\usepackage{amsmath}
\usepackage{a4}
\usepackage[english, francais]{babel}
\usepackage{euscript}
\usepackage{amsfonts}
\usepackage{amssymb}
\usepackage{mathrsfs}
\usepackage{xypic}
\usepackage[dvips]{graphics}
\usepackage{enumerate}
\usepackage{array}
\usepackage{theorem}
\usepackage{epic}
\usepackage{color}
\xyoption{all}
\usepackage[T1]{fontenc}
\usepackage{makeidx}
\usepackage{float}

\theoremstyle{plain}
\newtheorem{theorem}{Théorème}[section]
\newtheorem{proposition}[theorem]{Proposition}
\newtheorem{corollary}[theorem]{Corollaire}
\newtheorem{lemma}[theorem]{Lemme}

\theorembodyfont{\upshape}

\newtheorem{definition}[theorem]{Définition}

\newtheorem{remark}[theorem]{Remarque}

\newcommand{\qed}{\hfill $\Box$ \medskip}

\newcommand{\NN}{{\mathbb N}}
\newcommand{\ZZ}{{\mathbb Z}}
\newcommand{\QQ}{{\mathbb Q}}

\newcommand{\PP}{{\mathbb P}}
\renewcommand{\AA}{{\mathbb A}}

\newcommand{\GG}{\ensuremath{\mathbb G}}
\newcommand{\mmu}{\ensuremath{\boldsymbol \mu}}

\newcommand{\Spec}{{\operatorname{Spec}\kern 1pt}}
\newcommand{\R}{{\mathrm{R}}}
\renewcommand{\H}{{\mathrm{H}}}
\newcommand{\Aut}{{\mathrm{Aut}}}
\newcommand{\Isom}{{\mathrm{Isom}}}
\newcommand{\supp}{{\mathrm{supp}}}

\newcommand{\et}{{\textrm {\'et}}}

\newcommand{\p}{{\mathfrak p}}

\newcommand{\preuve}{\noindent{\it Preuve : }}

\newcommand{\limind}{\mathop{\mathop\mathrm {lim}\limits_{\xrightarrow{\hskip 0.5cm}}}}
\newcommand{\limproj}{\mathop{\mathop\mathrm {lim}\limits_{\longleftarrow}}}

\newcommand{\gen}{\mathfrak g}
\begin{document}

\title{Quelques calculs d'espaces $\R^i f_* G$ sur des courbes}

\author{Sylvain Maugeais}

\date{}

\maketitle

\begin{abstract}
Nous donnons quelques propriétés (nullité, représentabilité,
stratification) des espaces $\R^i f_* G$ pour $f\colon \mathcal U \to
S$ une courbe affine lisse possédant une compactification lisse et $G$
un groupe résoluble.   
\end{abstract}

\section{Introduction}

Soit $S$ un schéma, $\mathcal C \to S$ une courbe propre et lisse et
$\mathcal U \subset \mathcal C$ le complémentaire d'un diviseur de Cartier
horizontale. Notons $f\colon \mathcal U \to S$ le morphisme
induit. Pour tout groupe $G$ résoluble, nous calculons les faisceaux $\R^i
f_* G$ pour $i \in \{1, 2\}$.  Lorsque $f$ est propre, ces faisceaux
sont constructibles. Le cas le plus intéressant est donc le cas où
$\mathcal U$ est affine. D'autre part, si l'ordre de $G$ est
inversible sur $S$ alors le calcul est classique et peut se traiter
avec des méthodes habituelles (voir le chapitre $\S$
\ref{CalculModere} pour un exemple). Le premier cas intéressant est 
donc celui où $S$ est le spectre d'un corps algébriquement clos de
caractéristique $p$ divisant l'ordre de $G$, on a alors
$\R^2 f_* G=0$ et les $\R^1 f_* G$ peuvent être décris comme des
limites inductives et faisceaux représentables
(cf. \cite{Harbater}). Notre but est de généraliser ces résultats dans le
cas relatif. En particulier, nous montrons le théorème suivant, qui
constitue le coeur de la première partie de cet article 

\begin{theorem}
\label{MainTheorem}
Soit $S$ un schéma, $\mathcal C \to S$ une courbe propre et lisse et
$\mathcal U \subset \mathcal C$ le complémentaire d'un diviseur de Cartier
relatif tel que $\mathcal U \to S$ soit affine. Notons $f\colon
\mathcal U \to S$ le morphisme induit. Pour tout groupe $G$ on a
$\R^2 f_* G=0$.De plus si $G$ est résoluble il existe un
épimorphisme pour la topologie étale $\mathcal T \to \R^1 f_* G$ où
$\mathcal T$ est représentable par une limite inductive d'espaces
algébriques.  
\end{theorem}

En particulier, on trouve une généralisation du principe local-global
de \cite{KatzGabber} en la suite exacte
\begin{equation}
\label{PrincipeLocalGlobal}
0 \to \R^1 h_* G \to \R^1 f_* G \to h_* \R^1 j_* G \to \R^2 h_* G \to
0 
\end{equation}
Comme $h$ est propre, les termes $\R^i h_* G$ sont assez bien
maîtrisés. Il convient donc de mieux comprendre le terme $h_* \R^1 j_* G$.

Ce théorème se démontre par dévissage en se ramenant au cas
$\ZZ/p\ZZ$. L'un des problèmes fondamentaux est alors de
pouvoir traiter simultanément le cas $p$ inversible ou nul sur $S$ :
si $p$ est nul sur $S$ on montre que  $\R^2 f_* \ZZ/p\ZZ=0$ en
utilisant la suit exacte d'Artin-Schreier, mais si $p$ est inversible
c'est une conséquence du théorème de changement de base lisse. L'outil
principal utilisé ici est la suite exacte de déformation
d'Artin-Schreier à Kummer introduite par Oort, Sekiguchi et Suwa dans 
\cite{Oort_Sekiguchi_Suwa}. Le résultat s'obtient alors par une suite de
dévissage cohomologique.

Dans le cas où $S$ est le spectre d'un corps, le résultat
de\cite{KatzGabber} (qui correspond à la suite exacte
\ref{PrincipeLocalGlobal} dans le cas où la base est un corps) est
particulièrement utile pour construire des revêtements galoisiens de
courbes. Le cas de bases plus générales est toutefois beaucoup plus
compliqué. Pour l'étudier, nous introduisons une stratification
$\coprod_g (\R^1 f_* G)_{\gen=g}$ des espaces $\R^1 f_* G$ et nous
montrons que chaque strate est très proche de certains espaces de
modules de courbes équivariantes.  

Plus précisément étant donné un groupe fini $G$ et deux entiers $g$ et
$g'$, considérons l'espace des modules $\mathcal M_{g, g'}[G]$
classifiant les courbes propres et lisses de genre $g$ munie d'une
action de $G$ fidèle dans chaque fibre, dont le quotient est de genre
$g'$. Notons $n=|G|((2g-2)-|G|(2g'-2))$ et considérons l'espace
$\mathcal M^{[n]}_{g'}$ qui classifie les courbes propres et lisses de
genre $g'$ munie d'un diviseur de Cartier relatif de degré $n$. En
associant à chaque courbe équivariante son quotient muni du diviseur
de branchement, on obtient un morphisme  
$$\Phi\colon \mathcal M_{g, g'}[G] \to \mathcal M_{g'}^{[n]}.$$  

Notons $\mathfrak C$ la courbe universelle 
au-dessus de $\mathcal M_{g'}^{[n]}$ et $\mathfrak B$ le diviseur
universel. Notons de plus $\mathfrak f\colon \mathfrak C \setminus
|\mathfrak B| \to \mathcal M^{[n]}_{g'}$. 
On peut alors, à tout morphisme $S \to \mathcal M_{g, g'}[G]$,
associer un $G$-torseur au-dessus de $(\mathfrak C \setminus
|\mathfrak B|) \times_{\mathcal M_{g'}} S$ puis un morphisme de $S \to
[\R^1 \mathfrak f_* G]$, ce dernier espace étant un analogue champêtre
des $\R^1 f_* G$. On obtient ainsi un morphisme de champs $\mathcal
M_{g, g'}[G] \to [\R^1 \mathfrak f_* G]$ dont l'image est en fait dans
la strate $[\R^1 \mathfrak f_* G]_{\gen=g}$. 

En utilisant \cite{Alterations}, nous montrons qu'il existe une
altération $X \to [\R^1 \mathfrak f_* G]_{\gen=n}$ et un diagramme
commutatif  
$$\xymatrix{& X \ar[ld]\ar[d] \\ 
\mathcal M_{g, g'}[G] \ar[r] & [\R^1 \mathfrak f_* G]_{\gen=g}.}$$
Le cas d'un morphisme $f$ comme en début d'introduction s'obtient
alors par changement de base à partir de ce cas universel. 
D'autre part nous montrons que, quitte à remplacer $[\R^1
\mathfrak f_* G]$ par une des rigidifications données par le théorème
\ref{MainTheorem}, le morphisme $\mathcal M_{g, g'}[G] \to [\R^1
\mathfrak f_* G]_{\gen=g}$ est une immersion fermée.

Nous commençons donc au $\S$ \ref{EspacesTorseurs} par quelques
rappels sur les espaces de torseurs, en particulier des dévissages
cohomologiques et des applications de la suite exacte de déformation
d'Artin-Schreier à Kummer, et obtenons les premiers résultats sur les 
faisceaux $\R^i f_* G$ qui nous seront utiles par la suite. 

Dans $\S$ \ref{RepresentabiliteChamp} nous introduisons l'analogue
champêtre $[\R^i f_* G]$ des faisceaux $\R^i f_* G$ et montrons
l'existence d'une rigidification ind-représentable comme énoncée dans
le théorème \ref{MainTheorem}. Celle-ci permet alors de conclure sur
la nullité des $\R^2 f_* G$ dans le cas affine.

L'étude du morphisme $\Phi$ construit ci-dessus fait l'objet du
 $\S$ \ref{EtudeMorphisme} qui introduit proprement les
espaces de modules utilisé ci-dessus, la stratification par le genre
puis démontre les liens entre les espaces construits.

Finalement, nous utilisons les constructions de ce chapitre dans
l'appendice $\S$ \ref{CalculModere} pour retrouver simplement des
résultats connus sur l'espace des modules des courbes équivariantes
dans le cas modéré.

\noindent{\bf Remerciements :} L'auteur remercie Laurent Moret-Bailly,
Lorenzo Ramero, Michel Raynaud et  Matthieu Romagny et pour les
conversations qu'il a pu avoir avec chacun d'eux pendant la
préparation de cet article.

\section{Faisceaux de torseurs}
\label{EspacesTorseurs}

Le but de cette section est de donner quelques résultats sur les
torseurs et en particulier de généraliser des résultats de Harbater
(cf. \cite{Harbater}) et de Katz (cf. \cite{KatzGabber}). Plus
précisément, soit $C$ une courbe lisse sur un corps algébriquement
clos $k$, $U$ un sous-schéma ouvert de $C$ affine et $G$ un groupe
abstrait. L'un des buts de chacun de ces articles était de calculer le
groupe $\H^1_{\et}(U, G)$ et de donner ainsi un principe local-global
pour les revêtements. 

Nous souhaitons ici généraliser ce résultat dans le cas d'une courbe
relative $\mathcal C \to S$. Pour cela, l'un des instruments
essentiels est la suite exacte de déformation d'Artin-Schreier à
Kummer dont nous rappelons la théorie ci-après.

Dans toute la suite, lorsque les groupes ne sont pas commutatifs,  
les suites exactes de cohomologie sont à comprendre au sens non
abélien, i.e. suites exactes de faisceaux d'ensembles pointés.

\subsection{Rappels et notations}
\label{DefFaisceau}

Dans toute la suite de cet article, $S$ sera un schéma, $h\colon X \to
S$ un morphisme de schémas et $\mathcal G$ un schéma en groupes étale
sur $X$. 

Pour $i \le 2$, on peut alors définir des ensembles pointés
$\H^i_{\et}(\mathcal C, \mathcal G)$, cf. \cite{Giraud}. Nous noterons
$\R^i h_* \mathcal G$ le faisceau sur le gros site étale $S_{\et}$
associé au préfaisceau $$T \mapsto \H^i_{\et} (X \times_S T, \mathcal
G \times_S T).$$ 

\begin{proposition}
Supposons $h$ propre. Alors $\R^i h_* \mathcal G$ est représentable
par un espace algébrique pour $i \in \{1, 2\}$.
\end{proposition}

\preuve Pour $i=1$, le résultat s'obtient de la même manière que dans
le cas commutatif en se ramenant à des déformations de torseurs. 

Pour $i=2$, on sait que $\H^2_{\et}(X\times_S T, \mathcal G\times_S
T)$ est un espace principal homogène sous $\H^2_{\et}(X \times_S T,
Z(\mathcal G\times_S T))$ où $Z(\mathcal G\times_S T)$ désigne le
centre de $\mathcal G\times_S T$ (cf. \cite{Giraud}, Théorème
IV.3.3.3). Par suite, le faisceau $\R^2 h_* \mathcal G$ est un espace
principal homogène sous $\R^2 h_* Z(\mathcal G)$. On se ramène donc au
cas abélien, pour lequel le résultat est connu (cf. par exemple
\cite{Milne}). \qed 

Bien évidement, il est possible d'appliquer les opérations classiques
de la cohomologie. Nous nous contentons ici de donner les exemples qui
nous serons utiles (cf. \cite{Giraud}). 

Donnons nous une suite exacte de groupes étales sur $S$ 
$$1 \to \mathcal H \to \mathcal G \to \mathcal G/\mathcal H \to 1$$
avec $\mathcal H\subset Z(\mathcal G)$.
Alors il existe une suite exacte de faisceaux d'ensembles pointés 
\begin{multline}
\label{DevissageGroupe}
0 \to \R^0 h_* \mathcal H \to \R^0 h_* \mathcal G \to \R^0 h_*\mathcal
G/\mathcal H \to \R^1 h_* \mathcal H \to \R^1 h_* \mathcal G \to \R^1
h_*\mathcal G/\mathcal H \to \R^2 h_* \mathcal H
\end{multline}

Cette suite peut être utilisée pour l'étude des faisceaux de torseurs
sous des groupes résolubles afin de se ramener à l'étude des torseurs
sous des groupes cycliques. 

\medskip

Le dévissage sous forme de suite exacte comme ci-dessus est très
pratique en théorie mais difficile à manipuler en pratique. Dans les
faits, il est souvent plus utile d'avoir des caractérisations des
extensions galoisiennes d'extensions galoisiennes qui sont elles-mêmes
galoisiennes. Ceci est donné par la proposition suivante.

\begin{proposition}
\label{GalGalGal}
Soit $X$ un schéma, $P$ et $H$ des groupes finis, $Y$ un $P$-torseur
au-dessus de $X$ et $Z$ un $H$-torseur au-dessus de $Y$. Si $Z$ est 
un schéma connexe alors c'est un $\Aut_X(Z)$-torseur si et seulement
si $Z \in \H^1_{\et}(Y, H)^P$.  

De plus, on a une suite exacte $0 \to H \to \Aut_X(Z) \to P \to 0$
et si $Z$ et $Z'$ sont deux $H$-torseurs au-dessus de $Y$ qui sont
isomorphes, alors ils induisent des $\Aut_X(Z)$-torseurs isomorphes
au-dessus de $X$. 

En particulier, on a une applications injective 
$$\H^1_{\et}(Y, H)^P \to \coprod_{\{0 \to H \to G \to P \to 0\}/isom}
\H^1_{\et}(X, G)$$ 
\end{proposition}

On prendra garde au fait que $\H^1_{\et}(Y, H)^P$ n'est pas un groupe
en général, même si $H$ est abélien.

\preuve Pour tout $\sigma \in P$ notons $Z_\sigma$ le $H$-torseur
au-dessus de $Y$ obtenu par changement de base $\sigma\colon Y \to Y$
à partir de $Z$.  
%En général, le morphisme $Z_\sigma \to Z$ n'est pas un $Y$-morphisme

Supposons tout d'abord que $Z \to X$ est un $\Aut_X(Z)$-torseur. En
particulier, on a un morphisme surjectif $\Aut_X(Z) \to P$. N'importe
quel relèvement de $\sigma$ dans $\Aut_X(Z)$ fournit alors un
isomorphisme entre $Z_\sigma$ et $Z$.

Réciproquement, supposons que $Z \in \H^1_{\et}(Y, H)^P$. D'après
\cite{FundamentalGroup}, Lemma 2.8. il suffit de prouver que
$Z \times_X Z$ est totalement décomposé.  Comme $Y \to X$ est un
$P$-torseur, on a un isomorphisme 
$$Z \times_X Z \cong \prod_{\sigma \in P} Z_\sigma \times_Y Z.$$ 
Par suite, il suffit de montrer que $Z_\sigma \times_Y Z$ est un
produit de copies de $Z$. Or par hypothèse, $Z$ et $Z_\sigma$ sont
isomorphes en tant que $H$-torseurs donc $Z_\sigma \times_Y Z \cong Z
\times_Y Z \cong \prod_{\tau \in H} Z$. Ce qui fournit le résultat
annoncé. 

La dernière partie est claire.
% Si on a un isomorphisme de $H$-torseurs entre $Z$ et $Z'$, cet isomorphisme
% induit un isomorphisme entre $Aut_X(Z)$  et $\Aut_X(Z')$
\qed

\medskip

Donnons nous maintenant un morphisme propre $h\colon \mathcal C \to S$
ainsi qu'une immersion ouverte $j \colon \mathcal U \to \mathcal C$ et
notons $f\colon \mathcal U \to S$ le morphisme induit. On se fixe
également un schéma en groupes étale $\mathcal G$ sur $\mathcal U$. 

On a alors une suite exacte à $5$ termes 
\begin{equation}
\label{DevissageFaisceau}
0 \to \R^1 h_* (j_* \mathcal G) \to \R^1 f_* \mathcal G \to h_* \R^1
j_* \mathcal G \to \R^2 h_* (j_* \mathcal G) \to \R^2 f_* \mathcal G 
\end{equation}
obtenue en faisceautisant l'analogue non abélien de la suite exacte en
bas degré obtenue par la suite spectrale de Leray (cf. \cite{Giraud},
Chapitre V). 

Cette suite permet de décomposer $\R^1 f_* \mathcal G$ en un terme
local $h_* \R^1 j_*  \mathcal G$ et en termes globaux, cette suite est
l'un des outils fondamentaux utilisés par Harbater
(cf. \cite{Harbater}) et Katz (cf. \cite{KatzGabber}) lorsque $S$ est
le spectre d'un corps algébriquement clos. La suite exacte ci-dessus
doit être vue comme un principe local-global. Elle est toutefois
complètement différente de celle de Bertin et Mézard dans
\cite{BertinMezard} car elle ne dit rien sur les déformations
infinitésimales de courbes. 

À partir de maintenant, nous supposerons que $h$ est une courbe propre
et lisse et que $\mathcal U$ est le complémentaire dans $\mathcal C$
d'un diviseur de Cartier relatif. 

\subsection{La suite exacte d'Artin-Schreier à Kummer et quelques
  conséquences} 
\label{ArtinSchreier}

Nous allons maintenant rappeler quelques propriétés de la déformation
d'Artin-Schreier à Kummer (pour plus de détail, voir
\cite{Oort_Sekiguchi_Suwa}) et donnerons quelques conséquences, la 
plus importante étant la proposition \ref{Description-pcyclique}.

Dans toute la suite, nous fixons un nombre premier $p$ et un anneau de
valuation discrète $R_0$ de caractéristiques $(0, p)$ contenant une
racine $p$-ième de l'unité $\zeta$. Notons $\lambda=\zeta-1$. Il est
alors aisé de voir que $\lambda^{p-1}=up$ avec $u$ inversible dans
$R_0$. 

Pour tout $R_0$-schéma $X$ et tout $\mu \in R_0$ dans l'idéal maximal
de $R_0$, notons $\mathscr G^{(\mu)}_X=X \times_{R_0}\Spec R_0\left[x,
  \frac{1}{\mu x+1} \right]$. On 
peut munir $\mathscr G^{(\mu)}_X$ d'une structure de schéma en groupes
en utilisant la loi $a * b = \mu ab + a + b$. Sa partie en
caractéristique $0$ est alors isomorphe à $\GG_m$ et celle en
caractéristique $p$ à $\GG_a$. 

Notons $i_\mu$ l'immersion $X \times_{R_0} \Spec R_0/(\mu) \to X$. On
a une suite exacte pour la topologie \emph{fppf}  

\begin{equation}
\label{LienAvecGm}
0 \to \mathscr G^{(\mu)}_X \to \GG_m \to i_{\mu*} \GG_m \to 0
\end{equation}

On définit alors un morphisme de groupes $\psi:\mathscr
G^{(\lambda)}_X \to \mathscr G^{(\lambda^p)}_X$ par $x \mapsto
\frac{(\lambda x +1)^p -1 }{\lambda^p}$ qui est bien défini car
$\lambda^{p-1}=up$.   
On voit alors que $\psi$ est fidèlement plat et que son noyau est (non
canoniquement) isomorphe à $\ZZ/p\ZZ$. On a ainsi une suite exacte sur
le gros site étale
\begin{equation}
\label{ArtinShreierKummer}
0 \to \ZZ/p\ZZ \to \mathscr G^{(\lambda)}_X \to \mathscr
G^{(\lambda^p)}_X \to 0
\end{equation}
dont la fibre spéciale s'identifie à la suite exacte d'Artin-Schreier
et la fibre générique à la suite exacte de Kummer.

Appliquons cette théorie dans notre cas. Pour cela, considérons un
$R_0$-schéma $S$, une courbe propre et lisse $h\colon \mathcal C \to S$
et un ouvert $\mathcal U \subset S$ qui est le complémentaire d'un
diviseur de Cartier relatif. Notons $j\colon \mathcal U \to \mathcal
C$ l'inclusion et $i_\mu\colon\mathcal U\times_{R_0} \Spec
\R_0/(\mu) \to \mathcal U$.  Regardons tout d'abord la cohomologie de
$\GG_m$~:

\begin{lemma}
On a $\R^1 j_* \GG_m =  \R^1 j_* i_{\mu *} \GG_m = 0=(\R^2 j_* \GG_m
)_{tors}=(\R^2 j_* i_{\mu *}\GG_m)_{tors}$. 
\end{lemma}

\preuve La condition sur les termes de degré $1$ est en fait un
résultat de trivialité de faisceaux inversibles. Pour la montrer on peut
supposer que $S$ est local. Comme le résultat est local sur $\mathcal
C$, on peut remplacer ce dernier par un schéma local affine. Il s'agit
alors de montrer que tout faisceau inversible de $\H^1_{\et}(\mathcal
U, \GG_m)$ est trivial, c'est-à-dire que le torseur correspondant
possède une section. 

Notons $k$ le corps résiduel du point fermé de $S$. Comme $\mathcal
C_s$ est régulier, tout faisceau inversible de $\mathcal U_s$ peut se
prolonger en faisceau inversible sur $\mathcal C_s$. Or ces derniers
sont triviaux car $\mathcal C_s$ est local. Il est alors aisé de
passer au cas où $S$ est local en relevant une base.

Regardons maintenant les termes de degré $2$. D'après \cite{Milne},
Theorem III.3.9, on peut remplacer $S$ par le spectre d'un anneau
local complet $A$ à corps résiduel algébriquement clos d'idéal maximal 
$\mathfrak m$, $\mathcal C$ par $\Spec A[[x]]$ et $\mathcal U$ par
$\Spec A[[x]][\frac{1}{f}]$ où $f$ est une fonction régulière modulo
$\mathfrak m$. 

Si $A$ est un corps, alors la nullité provient de l'étude du groupe de
Brauer des corps locaux (cf. \cite{Corps_Locaux}).  Si $A$ est
artinien, la nullité de $\H^2_{et}(\Spec A[[x]][\frac{1}{f}], \GG_m)$
se démontre \emph{via} la théorie des déformations en utilisant le
fait que $\Spec A[[x]][\frac{1}{f}]$ est affine, ce qui implique que 
$\H^2_{et}(\Spec A[[x]][\frac{1}{f}], \GG_a)=0$.

Dans le cas général (complet), on a un morphisme canonique
$$\H^2_{et}(\Spec A[[x]][\frac{1}{f}], \GG_m) \to \limproj_n
\H^2_{et}(\Spec (A/\mathfrak m^n)[[x]][\frac{1}{f}], \GG_m)$$
dont il semble difficile de montrer l'injectivité. Toutefois, d'après
\cite{GabberBrauer}, le morphisme induit sur les parties de torsions
s'identifie à 
$$\mathrm{Br}(\Spec A[[x]][\frac{1}{f}]) \to \limproj_n
\mathrm{Br}(\Spec (A/\mathfrak m^n)[[x]][\frac{1}{f}]).$$
La démonstration de l'injectivité est alors semblable à
\cite{BrauerIII}, Lemme 3.3 en utilisant la trivialité de
$\H^1_{et}(\Spec A[[x]][\frac{1}{f}], \GG_m)$ démontrée ci-dessus.
\qed  

En appliquant $j_*$ à la suite exacte \eqref{LienAvecGm}, on trouve le
corollaire suivant.

\begin{corollary}
Avec les notations ci-dessus, on a 
$$\R^1 j_* \mathscr G^{(\mu)}_{\mathcal U} = (\R^2 j_* \mathscr
G^{(\mu)}_{\mathcal U})_{tors}=0.$$
\end{corollary}

\preuve Cela provient de la suite exacte \eqref{LienAvecGm} et de la
surjectivité du morphisme $j_* \GG_m \to j_* i_{\mu *} \GG_m$. \qed

Nous pouvons maintenant passer à la cohomologie de $\ZZ/p\ZZ$.

\begin{proposition}
\label{Description-pcyclique}
On a deux suites exactes pour la topologie étale 
\begin{equation}
\label{DevissageArtinSchreierKummerDegre1}
0 \to h_*\left((j_* \mathscr G^{(\lambda)}_U)/\mathscr
G^{(\lambda)}_{\mathcal C}\right) \to h_*\left((j_* \mathscr
G_{\mathcal U}^{(\lambda^p)})/\mathscr G^{(\lambda^p)}_{\mathcal C}
\right) \to h_*\R^1 j_* \ZZ/p\ZZ \to 0.
\end{equation}

\begin{equation}
\label{DevissageArtinSchreierKummerDegre2}
h_* \R^1 j_* \ZZ/p\ZZ \to \R^2 h_* \ZZ/p\ZZ \to \R^2 f_* \ZZ/p\ZZ \to 0
\end{equation}
\end{proposition}

\preuve En prenant la suite exacte longue de cohomologie associée à la suite
\eqref{ArtinShreierKummer} par application du foncteur $j_*$ on obtient une
suite exacte 
$$0 \to \ZZ/p\ZZ \to j_*\mathscr G_{\mathcal
  U}^{(\lambda)} \to j_*\mathscr G_{\mathcal U}^{(\lambda^p)} \to \R^1
j_*\ZZ/p\ZZ \to \R^1 j_* \mathscr G_{\mathcal U}^{(\lambda)}$$ 
dont le dernier terme est nul d'après le corollaire précédent. 

En particulier, on a un diagramme  commutatif à lignes exactes

$$\xymatrix{
0 \ar[r] & \ZZ/p\ZZ \ar[r]\ar[d] & \mathscr G_{\mathcal C}^{\lambda}
\ar[r]\ar[d] & \mathscr G_{\mathcal C}^{(\lambda^p)} \ar[r]\ar[d] & 0 \\ 
0 \ar[r] & \ZZ/p\ZZ \ar[r] & j_*\mathscr G_{\mathcal U}^{\lambda} \ar[r] &
j_*\mathscr G_{\mathcal U}^{(\lambda^p)} \ar[r] & \R^1 j_*
\ZZ/p\ZZ \ar[r] & 0. 
}$$
Le lemme du serpent fournit alors la suite exacte 

\begin{equation}
\label{SupportDim0}
0 \to (j_* \mathscr G_{\mathcal U}^{(\lambda)})/\mathscr
G^{(\lambda)}_{\mathcal C} \to (j_*
\mathscr G_{\mathcal U}^{(\lambda^p)})/ \mathscr
G^{(\lambda^p)}_{\mathcal C} \to \R^1 j_* \ZZ/p\ZZ \to
0.
\end{equation}

Les faisceaux présents dans cette suite étant à support fini sur $S$ (le
support est dans $\mathcal C \setminus \mathcal U$ qui est le support d'un
diviseur de Cartier relatif), on obtient la suite exacte
\eqref{DevissageArtinSchreierKummerDegre1} en appliquant le foncteur
$h_*$ qui est exact sur ces faisceaux.  
%% cf Scripte Milne, Corollary 8.4

\medskip

La suite exacte \eqref{DevissageArtinSchreierKummerDegre2} provient
de la suite spectrale de Leray
$$\R^\ell h_* \R^m j_* \ZZ/p\ZZ \Rightarrow \R^{\ell+m} f_*
\ZZ/p\ZZ$$ et d'un argument de dégénérescence \emph{via} le lemme
suivant.

\begin{lemma}
\label{NulliteH2Inclusion}
On a 
$$\R^3 h_* \ZZ/p\ZZ=\R^2 h_* \R^1 j_* \ZZ/p\ZZ=\R^1 h_* \R^1 j_*
\ZZ/p\ZZ = \R^2 j_* \ZZ/p\ZZ = 0.$$
\end{lemma}

\preuve On a $\R^3 h_* \ZZ/p\ZZ=0$ pour des raisons de dimension
cohomologique (cf. \cite{Milne}, Corollary VI.2.5). 

D'après la suite exacte \eqref{SupportDim0}, on voit que $\R^1 j_*
\ZZ/p\ZZ$ est à support de dimension relative nulle sur $S$ (car dans
$\mathcal C \setminus\mathcal U$), donc
$\R^\ell h_* \R^1 j_* \ZZ/p\ZZ=0$ pour tout $\ell > 0$.

La suite exacte longue associée à la suit exacte
\eqref{ArtinShreierKummer} nous donne
$$\R^1 j_*\mathscr G_{\mathcal U}^{(\lambda^p)} \to \R^2 j_*
\ZZ/p\ZZ \to \R^2 j_* \mathscr G_{\mathcal U}^{(\lambda)}.$$
Mais comme $\R^2 j_* (\ZZ/p\ZZ)_{\mathcal U}$ est de torsion, on a
en fait une suite exacte 
$$\R^1 j_*\mathscr G_{\mathcal U}^{(\lambda^p)} \to \R^2 j_*
\ZZ/p\ZZ \to (\R^2 j_* \mathscr G_{\mathcal U}^{(\lambda)})_{tors}.$$
Or on a vu que $\R^1 j_*\mathscr G_{\mathcal U}^{(\lambda^p)}=0$ et 
$(\R^2 j_* \mathscr G_{\mathcal U}^{(\lambda)})_{tors}=0$.
\qed

\section{Représentabilité et description}
\label{RepresentabiliteChamp}

Notre but est ici de généraliser aux champs algébriques les
constructions faites dans \ref{DefFaisceau} puis de montrer
l'existence d'une rigidification comme dans le théorème
\ref{MainTheorem}. Cette rigidification permettra alors de conclure
sur la nullité de $\R^2 f_* \ZZ/p\ZZ$. 

Fixons un champ algébrique $\mathcal S$ et un champ
algébrique en groupes représentable et étale $\mathcal G \to \mathcal
S$.  

\begin{definition}
Soit $f\colon\mathcal X \to \mathcal S$ un morphisme représentable de
champs algébriques. Notons $[\R^1 f_* \mathcal G]$ le groupo\"\i de
sur $\mathcal S$ dont les fibres au-dessus d'un schéma $T$ sont les
catégories dont  
\begin{enumerate}[-]
\item les objets sont les couples $(\alpha,\ \xi)$ avec $\alpha\colon
  T \to \mathcal S$ et $\xi \in \H^0(T, \R^1 f_{T *} \mathcal G_T)$ où
  $f_T$ et $\mathcal G_T$ sont obtenus de $f$ et $\mathcal G$ par
  changement de base ; 
\item les morphismes sont les automorphismes de $\alpha$ qui laissent
  fixe $\xi$. 
\end{enumerate}
\end{definition}

% L'hypothèse sur la representabilite de f est nécessaire pour pouvoir
% parler de \H^0(S, \R^1 f_{S *}) car si f n'est pas représentable,
% \mathcal U \times S n'est pas forcément un espace algébrique et ce
% n'est pas définis. 

Lorsque $\mathcal S$ est un espace algébrique, on retrouve le faisceau
$\R^1 f_* \mathcal G$ comme défini précédemment.

Le groupo\"\i de $[\R^1 f_* \mathcal G]$ doit être vu comme un espace
de module grossier ``partiel'' associé à $f_* [\mathcal U/\mathcal G]$
qui vit au-dessus de $\mathcal S$. En particulier, on a un morphisme
de $\mathcal S$-groupo\"\i de $$f_* [\mathcal U/\mathcal G] \to [\R^1
f_* \mathcal G].$$

Le groupo\"\i de $[\R^1 f_* \mathcal G]$ est en fait un champ (ceci
provient du fait qu'on considère le faisceau $\R^1 f_* \mathcal G$). 
%\preuve Soit $(\alpha, \xi)$ un objet de $[\R^1 f_* \mathcal G]$ et
%considérons $\Isom((\alpha, \xi))$. Soit $S' \to S$ un recouvrement
%étale de schémas et $h$ dans le noyau de $$\Isom((\alpha, \xi))(S')
%\to \Isom((\alpha, \xi))(S' \times_S S')$$ Comme $\mathcal S$ est un
%préchamp, $\Isom(\alpha)$ est un faisceau. L'élément $h$ définit donc
%un automorphisme de $\alpha$. Le résultat provient alors du fait que
%$\R^1 f_{S *} \mathcal G_S$ est un faisceau, donc il laisse fixe
%$\xi$ puisque c'est vrai localement. Ici on utilise juste
%l'injectivite dans la définition de faisceau \qed 
Par contre, même sous des hypothèses très restrictives, il n'est pas
algébrique.  Le problème fondamental est que le faisceau des
automorphismes d'un objet n'est pas représentable en général. 

Comme $[\R^1 f_* \mathcal G]$ n'est pas représentable en général, on
ne peut pas parler de séparation. Par contre, on peut vérifier le
critère valuatif de séparation. 

\begin{proposition}
\label{ChampSepare}
Le morphisme $[\R^1 f_* \mathcal G] \to \mathcal S$ vérifie le critère
valuatif de séparation. 
\end{proposition}

\preuve Donnons nous un anneau de valuation $R$ de corps de fractions
$K$ et un morphisme $\phi\colon\Spec K \to [\R^1 f_* \mathcal G]$ . 
Quitte à faire un changement de base étale, on peut supposer que
$\phi$ représente un $\mathcal G$-torseur $V \to \mathcal U
\times_{\mathcal S} \Spec K$. Par suite, si $V$ peut s'étendre en un
torseur $\mathcal V$ au-dessus de $\mathcal U\times_{\mathcal S} \Spec
R$, cela ne peut être que de manière unique car $\mathcal V$ est alors
le normalisé de $\mathcal U\times_{\mathcal S} \Spec R$ dans l'anneau
total des fractions de $V$. \qed

La suite exacte \eqref{DevissageFaisceau} suggère d'introduire un
analogue local de l'espace $[\R^1 f_* \mathcal G]$ :

\begin{definition}
Soient $\mathcal S$ un champ algébrique, $h\colon\mathcal C \to
\mathcal S$ une courbe propre et lisse et $j\colon\mathcal U \to
\mathcal C$ le complémentaire d'un diviseur de Cartier relatif
$\mathcal B$.  

On définit alors $[h_*\R^1 j_* \mathcal G]$ comme étant le groupoïde
dont les objets au-dessus d'un schéma $\mathcal S$-schéma $T$ sont les
éléments de $(\R^1 j_{T*} \mathcal G)(\mathcal U \times_{\mathcal S}
T)$ et les morphismes entre deux tels objets $\xi$ et $\xi'$ sont les
morphismes de $T \to \mathcal S$ qui envoient $\xi$ sur $\xi'$.
\end{definition}

Par localisation, on obtient un morphisme naturel 
$$[\R^1 f_* \mathcal G] \to [h_*\R^1 j_* \mathcal G]$$
et la suite exacte \eqref{DevissageFaisceau} permet d'en préciser
l'image et les fibres. Ce principe local-global permet en particulier
de montrer la proposition suivante dans laquelle on suppose $\mathcal
G = \ZZ/p\ZZ$. 

\begin{proposition}
\label{RepresentabiliteCasModere}
Soient $\mathcal S$ un champ algébrique, $h\colon\mathcal C \to
\mathcal S$ une courbe propre et lisse et $j\colon\mathcal U \to
\mathcal C$ le complémentaire d'un diviseur de Cartier relatif
$\mathcal B$.  Si l'entier $p$ est premier à toutes les
caractéristiques résiduelles de $\mathcal S$ alors $[\R^1 f_*
\ZZ/p\ZZ]$ est représentable. 

Si le champ $\mathcal S$ est le spectre d'un corps algébriquement clos
de caractéristique $p$, alors les points géométriques de $[\R^1 f_*
\ZZ/p\ZZ]$ sont les points d'une limite inductive de schémas :
l'espace des modules formels de Harbater.  
\end{proposition}

\preuve On se ramène aisément au cas où $\mathcal S$ est un espace
algébrique. Comme $h$ est propre, les faisceaux $\R^i h_* \ZZ/p\ZZ$
sont représentables. Utilisant la suite \eqref{DevissageFaisceau} on
se ramène à montrer que $[h_*\R^1 j_* \ZZ/p\ZZ]$ est représentable. 

Par suite la première partie est classique (voir la section
\ref{CalculModere}), la deuxième est démontrée par Harbater dans
\cite{Harbater}. \qed 

En particulier, même si $\mathcal S$ est un champ algébrique
quelconque, on peut parler de l'ensemble des points géométriques à
valeur dans $[\R^1 f_* \ZZ/p\ZZ]$ : il suffit de considérer les fibres
géométriques de $[\R^1 f_* \ZZ/p\ZZ] \to \mathcal S$ et l'ensemble de
leurs points géométriques, ce qui a un sens d'après la proposition
ci-dessus.  

Le problème fondamental pour la représentabilité de $[\R^1 f_*
\ZZ/p\ZZ]$ est que les déformations d'un torseur au-dessus d'un schéma
affine sont toutes isomorphes. Si $[\R^1 f_* \ZZ/p\ZZ]$ était
représentable, il serait étale sur $\mathcal S$, or on voit aisément
avec l'espace de modules formels d'Harbater que ce n'est pas le cas
(cf. \cite{Harbater}).

\subsection{Rigidification des $[\R^i f_* \ZZ/p\ZZ]$}

Nous allons maintenant introduire un espace proche de ce que devrait
être l'espace des déformations de $[h_*\R^1 j_* \ZZ/p\ZZ]$ si ce
dernier était représentable. Pour cela, l'outil fondamental est la
suite exacte \eqref{DevissageArtinSchreierKummerDegre1}. 

Dans toute la suite on se fixe un premier $p$ ainsi qu'une racine
primitive $p$-ième $\zeta$ de l'unité dans $\bar \QQ$. On notera de
plus $\lambda=\zeta-1$ (cf. \ref{ArtinSchreier}). 

\begin{definition}
Soient $\mathcal S$ un $\ZZ[\zeta]$-champ algébrique, $\mu \in
\ZZ[\zeta]$, $h\colon\mathcal C \to \mathcal S$ une courbe propre et
lisse et $j\colon\mathcal U \to \mathcal C$ le complémentaire d'un
diviseur de Cartier relatif. Nous noterons $\mathscr{T}_\mu(f)_{loc}$ le
$\mathcal S$-groupoïde dont les objets au-dessus d'un morphisme
$\alpha\colon T \to \mathcal S$ sont les éléments de $h_*\left((j_*
  \mathcal G_{\mathcal U}^{(\mu)})/\mathcal G^{(\mu)}_{\mathcal C} \right)(T)$ 
et les morphismes entre deux objets $\xi$ et $\xi'$ sont les
isomorphismes de $\alpha$ qui envoient $\xi$ sur $\xi'$.

En particulier, on a un morphisme $\phi\colon \mathscr{T}_{\lambda}(f)_{loc} \to
\mathscr{T}_{\lambda^p}(f)_{loc}$ 
(cf. \ref{ArtinSchreier}).
\end{definition}

On voit aisément que $\mathscr{T}_\mu(f)_{loc}$ est un champ. D'autre
part, il existe un morphisme canonique $$\mathscr{T}_{\lambda^p}(f)_{loc} \to
[h_*\R^1 j_* \ZZ/p\ZZ]$$ qui est un épimorphisme pour la topologie
étale d'après la suite exacte \eqref{DevissageArtinSchreierKummerDegre1}. 

Contrairement aux espaces précédemment introduit,
$\mathscr{T}_\mu(f)_{loc}$ est presque représentable. On effet, on a la
proposition suivante. 

\begin{proposition}
\label{RepresentabiliteLimite}
Le groupoïde $\mathscr{T}_\mu(f)_{loc}$ est une limite inductive de
champs algébriques de type finis et lisses sur $\mathcal S$.
\end{proposition}

\preuve Quitte à faire un changement de base lisse, on peut supposer
que $\mathcal S$ est en fait un schéma. Notons $\mathcal B$ un
diviseur de Cartier relatif tel que $\mathcal U = \mathcal C \setminus
|\mathcal B|$. Alors on voit qu'ensemblistement (les structures de
groupes n'étant pas les mêmes) on a $$h_*\left((j_* \mathcal
  G_{\mathcal U}^{(\lambda^p)})/\mathcal G^{(\lambda^p)}_{\mathcal C}
\right)= \limind_n h_*\left(\O_{\mathcal C}(n\mathcal B)/\O_{\mathcal
    C} \right)$$ et ce dernier est représentable car $\mathcal B$ est
fini et plat sur $\mathcal S$. 

Finalement la lissité provient de la platitude de $\O_C(n\mathcal
B)/\O_{\mathcal C}$ et du fait que le faisceau de ses sections est
représentable par  $\Spec \left(\mathrm{Sym} h_*\left(\O_C(n\mathcal
    B)/\O_{\mathcal C}\right)^\vee\right)$. \qed 

On peut donc voir $\mathscr{T}_{\lambda^p}(f)_{loc}$ comme une
rigidification de $[h_*\R^1 j_* \ZZ/p\ZZ]$. Nous allons maintenant
introduire une rigidification de $[\R^1 f_* \ZZ/p\ZZ]$ en utilisant
$\mathscr{T}_{\lambda^p}(f)_{loc}$ et la suite exacte
\eqref{DevissageFaisceau}. 

\begin{definition}
Soient $\mathcal S$ un $\ZZ[\zeta]$-champ algébrique, $h\colon
\mathcal C \to \mathcal S$ une courbe propre et lisse et $j\colon
\mathcal U \to \mathcal C$ le complémentaire d'un diviseur de Cartier
relatif. On note $$\mathscr{T}_{\lambda^p}(f) = [\R^1 f_* \ZZ/p\ZZ]
\times_{[h_*\R^1 j_* \ZZ/p\ZZ]} \mathscr{T}_{\lambda^p}(f)_{loc}.$$ 
\end{definition}

Si $\mathcal S$ est un espace algébrique, on a une suite exacte 
$$ 0 \to \R^1 h_* \ZZ/p\ZZ \to \mathscr{T}_{\lambda^p}(f) \to
\mathscr{T}_{\lambda^p}(f)_{loc} \to \R^2 h_* \ZZ/p\ZZ.$$ 
D'autre part, comme $h$ est propre, $\R^i h_* \ZZ/p\ZZ$ est
représentable pour tout $i \ge 0$. La proposition précédente montre
alors que $\mathscr{T}_{\lambda^p}(f)$ est un champ et est limite inductive de
champs algébriques. 

\begin{proposition}
\label{RelevementZariski}
Le morphisme canonique $\mathscr{T}_{\lambda^p}(f) \to [\R^1 f_*
\ZZ/p\ZZ]$ est un épimorphisme pour la topologie étale et son noyau
est $\mathscr T_{\lambda}(f)_{loc}$. 
\end{proposition}

\preuve Par descente, on se ramène au cas où $\mathcal S$ est un
espace algébrique. Les champs considérés sont alors des faisceaux. Le
résultat découle du diagramme commutatif à lignes exactes 
$$\xymatrix{
0 \ar[r] & \R^1 h_* \ZZ/p\ZZ \ar[r] \ar@{=}[d] &
\mathscr{T}_{\lambda^p}(f) \ar[r] \ar[d] &
\mathscr{T}_{\lambda^p}(f)_{loc} \ar[r] \ar[d] & \R^2 h_* 
\ZZ/p\ZZ \ar@{=}[d]\\ 
 0 \ar[r] & \R^1 h_* \ZZ/p\ZZ \ar[r]  & \R^1 f_* \ZZ/p\ZZ  \ar[r] &
 h_* \R^1 j_* \ZZ/p\ZZ  \ar[r]  & \R^2 h_* \ZZ/p\ZZ 
}$$
et du fait que le morphisme $\mathscr{T}_{\lambda^p}(f)_{loc} \to h_*
\R^1 j_* \ZZ/p\ZZ$ est un épimorphisme pour la topologie étale (voir
la suite exacte \eqref{DevissageArtinSchreierKummerDegre1}). \qed 

\begin{remark}
\label{RigidificationExplicite}
Dire que le morphisme $\mathscr{T}_{\lambda^p}(f) \to [\R^1 f_*
\ZZ/p\ZZ]$  est un épimorphisme revient à dire que si on se donne un 
$\ZZ/p\ZZ$-torseur on peut l'écrire localement sous la forme
$\frac{(\lambda u + 1)^p - 1}{\lambda^p}=t$. 

D'autre part, deux morphismes $T \to \mathscr{T}_p(f)$ sont égaux si
et seulement s'ils correspondent localement à des torseurs
$\frac{(\lambda u + 1)^p - 1}{\lambda^p}=t_1$ et $\frac{(\lambda u +
  1)^p - 1}{\lambda^p}=t_2$ avec $t_1-t_2$ n'ayant pas de pôles
(i.e. on ne prend pas de quotients supplémentaires). 
\end{remark}

Nous pouvons maintenant donner une première application de ces
résultats. 

\begin{theorem}
Supposons que $f$ est affine et que $\mathcal S$ est un schéma. Alors
$\R^2 f_* \ZZ/p\ZZ = 0$. 
\end{theorem}

\preuve D'après la suite exacte
\eqref{DevissageArtinSchreierKummerDegre2}, il suffit de montrer que
le morphisme $f_* \R^1 j_* \ZZ/p\ZZ \to \R^2 h_* \ZZ/p\ZZ$ est
un épimorphisme.

On a un diagramme commutatif dont la première ligne est exacte 
$$\xymatrix{
f_* \R^1 j_* \ZZ/p\ZZ \ar[r] & \R^2 h_* \ZZ/p\ZZ \ar[r] & \R^2 f_* \ZZ/p\ZZ\\
\mathscr T_{\lambda^p}(f)_{loc} \ar[u] \ar[ru]
}$$
D'autre part, le morphisme $f_* \R^1 j_* \ZZ/p\ZZ \to \R^2 h_*
\ZZ/p\ZZ$ est surjectif au niveau des points géométriques car $\R^2
f_* \ZZ/p\ZZ=0$ en chaque point (car $f$ est affine). 
Par suite, il en est de même du morphisme $\mathscr T_{\lambda^p}(f)_{loc} \to
\R^2 h_* \ZZ/p\ZZ$. Comme $\mathscr T_{\lambda^p}(f)_{loc}$ est
ind-représentable par des schémas lisses (cf. proposition
\ref{RepresentabiliteLimite}) et que $\R^2 h_* \ZZ/p\ZZ$ est
constructible, on trouve que $\mathscr T_{\lambda^p}(f)_{loc} \to \R^2 h_*
\ZZ/p\ZZ$ est un épimorphisme de faisceaux étales.
%% pour montrer ca, on se ramène au cas où $\R^2 h_* \ZZ/p\ZZ$ est de
%% la forme $i_! G$ avec $i:U \to X$ partie localement fermée et G
%% constant (cf. \cite{SGAIV}, exposé IX prop 2.5) 
%% comme surjectif en un point et lisse, on relève cette section dans
%% $\mathscr T_{p}(f)_{loc}$ et elle a la bonne image car elle a la
%% bonne image à la fibre spéciale et $\R^2 h_* \ZZ/p\ZZ$ loc constant !!!
\qed

\begin{corollary}
Supposons que $f$ est affine et que $\mathcal S$ est un schéma. Alors
pour tout groupe $G$ on a $\R^2 f_* G = 0$. 
\end{corollary}

\preuve D'après \cite{Giraud}, Théorème IV.3.3.3, il suffit de montrer
que $\R^2 f_* Z(G)=0$. On peut donc supposer que $G$ est abélien. Le
résultat découle alors aisément du cas $\ZZ/p\ZZ$ en utilisant la
suite \eqref{DevissageGroupe}. \qed

\begin{corollary}
Soit $G$ un groupe et $H \subset Z(G)$ un sous-groupe. Alors on a des
suites exactes   
$$0 \to h_* \R^1 j_* H \to h_* \R^1 j_* G \to h_* \R^1 j_* G/H \to 0$$
$$0 \to \R^1 f_* H \to \R^1 f_* G \to \R^1 f_* G/H \to 0$$
\end{corollary}

\preuve Cela découle de la suite exacte longue de cohomologie associée
à la suite exacte 
$$0 \to H \to G \to G/H \to 0$$ 
et du lemme \ref{NulliteH2Inclusion} qui assure que $\R^2 j_* G/H = 0$.
\qed

\begin{corollary}
\label{ExistenceRigidification}
Soit $G$ un groupe résoluble. Alors il existe une suite exacte (pour
la topologie étale) de champs en ensembles pointés
$$0 \to \mathscr T'_G(f) \to \mathscr T_G(f) \to [R^1 f_* G]\to 0$$
où $\mathscr T_G(f)$ et $\mathscr T'_G(f)$ sont des limites inductives
de champs algébriques de types finis lisses sur $\mathcal S$.
\end{corollary}

\preuve Provient du cas $\ZZ/p\ZZ$ et du corollaire ci-dessus par une
récurrence immédiate. \qed 

\subsection{Un ersatz de produit}

Nous aurons besoin dans le chapitre suivant de considérer des produits
de la forme 
$$\mathscr{T}_G(f) \times_{[\R^1 f_* G]} T$$ 
pour un morphisme $T \to [\R^1 f_* G]$ donné. Le problème est
que de tels produits ne sont pas représentables en général. Nous
allons donc introduire maintenant une alternative utilisant les
propriétés du morphisme $\mathscr{T}_G(f) \to [\R^1 f_* G]$ démontrée
dans le corollaire \ref{ExistenceRigidification}.

Fixons un champ algébrique $T$, un morphisme de champs $T \to [\R^1
f_* G]$ et supposons qu'il existe un relèvement $\phi\colon T
\to \mathscr{T}_G(f)$. On définit alors le champ algébrique
$\mathscr{T}_G(f) \tilde \times_{[\R^1 f_* G]} T$ comme étant
l'image schématique du morphisme 
$$\begin{array}{ccc}
\mathscr{T}'_G(f) \times T & \to & \mathscr{T}_G(f) \times T \\
(a, b) & \mapsto & (a+\phi(b), b).
\end{array}$$

On voit que $\mathscr{T}_G(f) \tilde \times_{[\R^1 f_* G]} T$
ainsi que l'inclusion $\mathscr{T}_G(f) \tilde \times_{[\R^1 f_*
  G]} T \to \mathscr{T}_G(f) \times T$ sont indépendants du
relèvement $\phi$ car deux relèvements diffèrent d'un élément de
$\mathscr{T}'_G(f)$. Cette construction est donc canonique. D'autre
part, elle commute aux changements de bases plats mais pas aux
changements de bases quelconques \emph{a priori} (car c'est le cas de
l'image schématique). 

On peut alors construire un champ algébrique $\mathscr{T}_G(f) \tilde
\times_{[\R^1 G]} T$ dans le cas général par recollement
car tout morphisme $T \to [\R^1 f_* G]$ peut être relevé
localement pour la topologie étale. On obtient de plus une immersion
fermée $\mathscr{T}_G(f) \tilde \times_{[\R^1 f_* G]} T \to
\mathscr{T}_G(f) \times T$.  

\medskip
Remarquons que si $[\R^1 f_* G]$ est représentable,
alors il est séparé d'après la proposition \ref{ChampSepare} et donc
le morphisme $\mathscr{T}_G(f) \times_{[\R^1 f_* G]} T \to
\mathscr{T}_G(f)  \times T$ est une immersion fermée. On a donc une
immersion fermée $\mathscr{T}_G(f) \tilde \times_{[\R^1 f_* G]}
T \to \mathscr{T}_G(f) \times_{[\R^1 f_* G]} T$ qui est une
bijection. Seule la structure de sous-schéma fermé de
$\mathscr{T}_G(f) \times T$ est éventuellement différente. 

On a finalement la propriété suivante qui est une conséquence directe
du fait que $\mathscr{T}_G(f) \to [\R^1 f_* G]$ est un
épimorphisme pour la topologie étale. 

\begin{proposition}
Le morphisme canonique $\mathscr{T}_G(f) \tilde \times_{[\R^1 f_*
  G]} T \to T$ est un épimorphisme pour la topologie étale.
\end{proposition}

\section{Extensions de torseurs et courbes équivaraintes}
\label{EtudeMorphisme}

Le but de ce chapitre est de comparer les espaces de torseurs
introduit précédemment à des espaces de modules de courbes
équivariantes. 

Dans un premier temps, nous donnons quelques résultats préliminaires,
en particulier une propriété d'extensions de morphismes de torseurs en
morphismes entre courbes lisses. Nous passons ensuite à quelques
propriétés du genre des courbe affines, ce qui nous permet de définir
une stratification des espaces $[\R^1 f_* G]$. nous comparons ensuite
les strates obtenues à des espaces de modules de courbes.

\subsection{Préliminaires}

Nous rassemblons ici quelques résultats techniques qui jouerons des
rôles fondamentaux dans la suite.

\begin{theorem}
\label{RelevementImpliqueIso}
Soient $Y$ un schéma localement noethérien et $f \colon X \to
Y$ un morphisme localement de type fini. Alors $f$ est une immersion fermée
si et seulement si les conditions suivantes sont vérifiées 
\begin{enumerate}[i)]
\item pour tout corps algébriquement clos $k$, l'application induite
  $X(\Spec k) \to Y(\Spec k)$ est injective ; 
\item le morphisme $f$ vérifie le critère valuatif de propreté ;
\item pour tout corps algébriquement clos $k$, l'application induite
  $X(\Spec k[\epsilon]/(\epsilon^2)) \to Y(\Spec
  k[\epsilon]/(\epsilon^2))$ est injective 
\end{enumerate}
\end{theorem}

\preuve Le sens direct est évident. Réciproquement, d'après la
propriété $ii)$, le morphisme $f$ est propre. La propriété $i)$ impose
alors que $f$ induit une bijection entre $X$ et un sous-schéma fermé
de $Y$. Par suite, la proposition $iii)$ montre que $f$ est une
immersion fermée.  \qed  

\medskip

Nous souhaitons regarder les courbes affines sur des bases
relativement quelconques qui possèdent des compactifications
lisses. L'outil principal pour cela est le théorème d'altérations de
de Jong.

\begin{theorem}
\label{ExistenceCompactification}
Soit $S$ un schéma, $\mathcal V \to S$ une courbe lisse non
nécessairement propre. Supposons que 
\begin{enumerate}[i)]
\item les fibres géométriques de $\mathcal V \to S$ sont connexes ;
\item il existe un entier $g$ tel que pour tout point géométrique
  $\bar s \to S$, l'unique complétion lisse de $\mathcal V_{\bar s}$
  est de genre $g$.
\end{enumerate}
Alors il existe une altération $S' \to S$ et une courbe propre
et lisse $\mathcal D \to S'$ contenant $\mathcal V \times_S S'$.
\end{theorem}

\preuve C'est un corollaire de \cite{Alterations} qui permet de
trouver une altération et une compactification semi-stable. Les
hypothèses faites permettent alors de voir que la compactification
semi-stable est en fait lisse. \qed

Dans le cas général, la compactification n'a aucune raison d'être
unique : si c'est le cas pour un corps, ce n'est déjà plus vrai pour
des anneaux artiniens. Toutefois, si $S'$ est normal, il est alors
aisé de montrer que $\mathcal D$ est unique à unique isomorphisme 
près. Le cas réduit semble plus difficile, le problème principal
venant du fait que pour une courbe propre et lisse $\mathcal D \to S$,
le schéma des automorphismes $\Aut_S(\mathcal D)$ est parfois ramifié
sur $S$. Par contre, ce ne sera plus le cas si $\mathcal D \to S$ est
de genre $\ge 2$ (cf. \cite{Deligne_Mumford}).

Il est possible de prouver l'unicité lorsqu'un morphisme fini vers une
courbe lisse est fixé. C'est l'objet du corollaire
\ref{ProlongementMorphisme}. Avant cela, nous avons besoin d'un lemme
précisant la géométrie du schéma des automorphismes.

\begin{lemma}
Soit $S$ un schéma, $\mathcal C \to S$ une courbe propre et lisse et
pour $i \in \{1, 2\}$, $\pi_i\colon \mathcal D_i \to \mathcal C$ des
$S$-morphismes entre courbes lisses, séparables dans chaque fibre
au-dessus de $S$. Alors le schéma des isomorphismes  
$\Isom_{\mathcal C}(\mathcal D_1, \mathcal D_2)\subset
\Isom_S(\mathcal D_1, \mathcal D_2)$ classifiant les $S$-isomorphismes
$\mathcal D_1 \to \mathcal D_2$ commutants avec $\pi_1$ et $\pi_2$, est
fini et non ramifié sur $S$. 
\end{lemma}

\preuve Tout d'abord, on voit que $\Isom_{\mathcal C}(\mathcal D_1,
\mathcal D_2)\subset \Isom_S(\mathcal D_1, \mathcal D_2)$ est
représentable sur $S$ car les courbes sont propres 
sur $S$. Montrons qu'il est quasi-fini sur $S$. Pour cela, on peut
supposer que $S$ est le spectre d'un corps. Les corps de fractions de
$\mathcal D_1$ et $\mathcal D_2$ sont finis sur le corps de fractions
de $\mathcal C$, il n'existe donc qu'un nombre fini d'isomorphismes. 

Par suite, il suffit de montrer qu'il est propre et non ramifié. 

Pour ce qui est de la propreté, il suffit de vérifier le critère
valuatif car on sait déjà qu'il est de présentation finie sur $S$. On
est donc ramené au cas où $S$ est le spectre d'un anneau de
valuation. Par suite, si on a un isomorphisme à la fibre générique, on
peut l'étendre car $\mathcal D_1$ et $\mathcal D_2$ sont les
normalisations de $\mathcal C$ dans leur corps de fractions.  

Passons maintenant à la non ramification. Pour cela, on peut supposer
que $S$ est le spectre d'un anneau artinien $A$ de corps résiduel
algébriquement clos, que $\mathcal D_1 = \mathcal D_2$ et qu'on a
$\phi\in\Isom_{\mathcal C}(\mathcal D_1, \mathcal D_1)$ qui est
l'identité modulo un élément $\epsilon \in A$ annulé par l'idéal
maximal de $A$. 

Notons $D_1$ la fibre spéciale de $\mathcal D_1$. Utilisant la théorie
des déformations des courbes, on voit que $\phi$ définit un élément de
$\chi \in \H^0(D_1, \Omega_{D_1}^\vee)$. Soit $\p \in D_1$ un point en
lequel le morphisme $\pi_1$ est étale. Comme $\phi$ est un $\mathcal
C$-morphisme, on trouve que $\chi$ possède un zéro en $\p$. Le
morphisme $\pi$ étant génériquement étale, le champ de vecteur $\chi$
s'annule sur un ensemble dense. Il s'ensuit que $\chi=0$ et donc $\phi=Id$. \qed

\begin{corollary}
\label{ProlongementMorphisme}
Soit $S$ un schéma noethérien \emph{réduit}, $\mathcal C \to S$, $\mathcal
D_1 \to S$ et $\mathcal D_2 \to S$ trois courbes propres et lisses, 
$\pi_1\colon\mathcal D_1 \to \mathcal C$ et $\pi_2\colon \mathcal D_2
\to \mathcal C$ deux revêtements finis tels que leur lieu de
branchement soit contenu dans le support d'un diviseur de Cartier
relatif $B \subset \mathcal C$. En particulier, on suppose que $\pi_1$
et $\pi_2$ sont génériquement séparables dans chaque fibre au-dessus
de $S$. Notons $\mathcal V_i:=\pi_i^{-1}(\mathcal C \setminus B)$ et
donnons nous un $\mathcal C$-isomorphisme $\phi\colon \mathcal V_1
\to \mathcal V_2$ au-dessus de $\mathcal C \setminus B$.  
Alors $\phi$ s'étend en un \emph{unique} $\mathcal C$-isomorphisme
$\mathcal D_1 \to \mathcal D_2$. 
\end{corollary}

\preuve Il s'agit de construire une section $S \to \Isom_{\mathcal
  C}(\mathcal D_1, \mathcal D_2)$ induisant $\phi$. 

Supposons dans un premier temps que $S$ est le spectre d'un corps ou
bien le spectre d'un anneau de valuation, ou plus généralement si $S$
est normal. Alors le résultat est trivialement vrai car $\phi$ induit
un isomorphisme sur les normalisations de $\mathcal C$ dans les corps
de fractions de $\mathcal V_1$ et de $\mathcal V_2$, comme celles-ci
sont $\mathcal D_1$ et $\mathcal D_2$, on obtient le résultat.

Notons $S_{gen}$ la réunion disjointes des points génériques de
$S$. Le schéma $S_{gen}$ est alors la réunion disjointe de spectres de
corps. Par suite, il existe une unique section $S_{gen} \to
\Isom_{\mathcal C}(\mathcal D_1, \mathcal D_2)$ induisant $\phi$
au-dessus de $S_{gen}$. Notons $Z$ l'adhérence schématique de cette
section dans $\Isom_{\mathcal C}(\mathcal D_1, \mathcal D_2)$ (la
structure de schéma étant la structure réduite). Il s'agit alors de
voir que $Z \to S$ est un isomorphisme. 

Comme $\Isom_{\mathcal C}(\mathcal D_1, \mathcal D_2) \to S$ est fini
et non ramifié, il en est de même de $Z \to S$. 

Soit $s \in S$ et $R$ un anneau de valuation dominant $\O_{S,
  s}$. Comme vu précédemment, on peut étendre l'isomorphisme à la
fibre générique (donné par $Z \times_{S} \Spec R$) par
normalisation. En particulier, on voit qu'il y a un unique point de
$Z$ au-dessus de $s$, et celui-ci est de même corps résiduel que $s$
car la section de $\Isom_{\mathcal C}(\mathcal D_1, \mathcal D_2)$
correspondante est obtenue par normalisation. 

Par suite, l'application ensembliste $Z \to S$ est injective. Comme ce
morphisme est non ramifié et propre, c'est une immersion fermée. 
Comme il est dominant et que $S$ est réduit, c'est un
isomorphisme. \qed 

En suivant la même preuve, on montre l'unicité de la compactification
d'une courbe affine $\mathcal D$ lorsque celle-ci est de genre $\ge 2$
dans le théorème \ref{ExistenceCompactification}.

\begin{corollary}
\label{ExtensionTorseur}
Soit $S$ un schéma, $\mathcal C \to S$ un morphisme propre et lisse,
$\mathcal U \subset \mathcal C$ le complémentaire d'un diviseur de
Cartier relatif et $\mathcal V \to \mathcal U$ un morphisme
étale. Supposons que 
\begin{enumerate}[i)]
\item les fibres géométriques de $\mathcal V \to S$ sont lisses et
  connexes ;
\item il existe un entier $g$ tel que pour tout point géométrique
  $\bar s \to S$, l'unique complétion lisse de $\mathcal V_{\bar s}$
  est de genre $g$.
\end{enumerate}
Alors il existe une altération $S' \to S$ et une \emph{unique} courbe
propre et lisse $\mathcal D \to S'$ contenant $\mathcal V \times_S S'$
telle qu'on a un diagramme commutatif
$$\xymatrix{
\mathcal V \times_S S' \ar[r] \ar[d] & \mathcal D \ar[d] \\
\mathcal U \times_S S' \ar[r] & \mathcal C \times_S S'
}$$
\end{corollary}

Il est possible d'étendre de résultat dans le cas modéré quand on
suppose que $\mathcal C \setminus \mathcal U$ est le support d'un
diviseur de Cartier relatif étale sur la base (cf. \S
\ref{CalculModere}). Le cas sauvage reste par contre difficile à
maîtriser. À titre d'exemple, considérons un corps $k$ non parfait de
caractéristique $p$, $a \in k \setminus k^p$ et $V$ le revêtement de
$\AA^1_k \setminus \{0\}$ donné par $y^p-y=\frac{a}{x^p}$. La
normalisation de $\PP^1$ dans le corps des fractions de $V$ est donnée
par $u^p-x^{p-1}u=a$ qui n'est pas lisse au-dessus de $x=0$.
Toutefois, il va exister une compactification lisse sur $k(\sqrt[p]{a})$.

Le corollaire ci-dessus est trivialement faux si $S$ n'est pas
réduit. En effet, les déformations infinitésimales de torseurs
au-dessus d'un schéma affine sont toujours isomorphes alors qu'il
existe en général des déformations infinitésimales de courbes propres
qui ne sont pas isomorphes. 

C'est ce phénomène qui explique que l'exemple ci-dessus ne peut pas se
redescendre à $k$ : l'anneau $k(\sqrt[p]{a}) \otimes_k
k(\sqrt[p]{a})$ n'étant pas réduit, on ne peut pas utiliser l'unicité
de la courbe pour redescendre la courbe lisse par descente finie plate.

\subsection{Stratification de l'espace des torseurs}

Soient $S$ un schéma et $\mathcal V \to S$ une courbe lisse, non
nécéssairement propre.
Pour tout point $s \in S$, nous noterons $\gen_{\mathcal V}(s)$ le
genre de l'unique compactification lisse de $\mathcal V_{\bar s}$, où
$\bar s$ est un point géométrique au-dessus de $s$. Cet entier ne
dépend en fait que de $s$ et commute à tout changement de base $S' \to
S$.

\begin{proposition}
Supposons que $S$ est un schéma noethérien, alors l'application
$\gen_{\mathcal V}\colon S \to \NN$ est constructible. 
\end{proposition}

\preuve Comme l'image d'un ensemble constructible par un morphisme de
schémas est constructible, on peut supposer que $S$ est normal et
intègre de point générique $\eta$. D'autre part, quitte à faire une
extension finie de $S$, on peut supposer que $\mathcal V_\eta$ possède
une compactification lisse (en général, cette extension est
nécessaire). Choisissons une compactification $\mathcal D \to S$ de
$\mathcal V$, qui existe d'après Nagata. En particulier, on peut la
choisir normale, ce qui assure que la fibre générique de $\mathcal D
\to S$ est lisse car on a supposé que $\mathcal V_\eta$ possède une
compactification lisse, celle-ci est automatiquement $\mathcal D_\eta$.

D'après \cite{EGA} Proposition IV.17.7.11, le lieu $S'$ des points $s
\in S$ où $\mathcal D_s \to \Spec k(s)$ est lisse est
constructible, et il est aisé de voir que $\gen_{\mathcal V}$ est
localement constant sur cet ensemble car c'est le cas de du genre de
$\mathcal D$. On peut alors répéter cette construction sur $S  
\setminus S'$ qui est de codimension $> 0$. Comme $S$ est noethérien,
on obtient le résultat après un nombre fini d'étapes.\qed 

En particulier pour tout entier $g \in \NN$ nous pouvons
considérer le sous-schéma $S_{\gen_{\mathcal V}=g}$ composé des points $s$ de
$S$ pour lesquels $\gen_{\mathcal V}(s)=g$. En effet, comme $S$ est noethérien
et $S_{\gen_{\mathcal V}=g}$ est constructible, ce dernier est une réunion
d'ensembles localement fermés. On peut donc le munir d'une structure
de schéma en considérant sa structure réduite. On prendra garde au
fait que ce schéma ne commute pas au changement de base quelconque car
la structure réduite n'a ici rien de canonique.

\medskip

Dans les faits, nous devrons utiliser cette construction pour un
groupoïde qui est limite inductive de champs algébriques
noethériens. Le fait que ce soit un groupoïde ne pose pas de problème
car l'application $\gen_{\mathcal V}$ commute aux changements de base. On
généralise alors aisément la construction ci-dessus au cas de la
limite inductive. 

\medskip

Nous donnons maintenant quelques propriétés de la différente. 

\begin{proposition}
\label{DifferenteDecroit}
Reprenons les notations précédentes et supposons $S$
noethérien. Soient $\eta$ un point de $S$ et $s$ une spécialisation de
$\eta$. Alors $\gen_{\mathcal V}(s) \le \gen_{\mathcal
  V}(\eta)$. 
\end{proposition}

\preuve On peut supposer que $S$ est le spectre d'un anneau de
valuation discrète, les points $\eta$ et $s$ étant ses points
génériques et fermés. Quitte à faire une extension finie de $S$, on
peut de plus supposer que les compactifications normales de chaque
fibre sont lisses. Par suite, on peut choisir une compactification
$\mathcal D$ de $\mathcal U$ telle que cette dernière est dense dans
chaque fibre. Si on suppose de plus $\mathcal D$ normale alors il est
Cohen-Macaulay (d'après le critère de Serre), et sa fibre spéciale est
réduite. Il est alors aisé de montrer que la caractéristique d'Euler de
la normalisée de la fibre spéciale est plus petite que la
caractéristique d'Euler de la fibre spéciale. Cette dernière étant
égale à la caractéristique d'Euler de la fibre générique par
platitude, on obtient le résultat souhaité. \qed

\subsection{Un morphisme entre espaces de modules de courbes}

Soient $g$ un entier et $G$ un groupe abstrait fini. Nous noterons
$\mathcal M_g[G]$ l'espace des modules des courbes propres et lisses
de genre $g$ munies d'une action de $G$ fidèle dans chaque fibre. Cet
espace a été introduit par Tufféry dans \cite{DeformationEquivariante}
et a été étudié indépendemment pas plusieurs auteurs. 

Soient $\mathcal D \to S$ une telle courbe et $H$ un sous-groupe de
$G$. Dans \cite{Samuel}, Samuel montre que le quotient $\mathcal D \to
\mathcal D/H$ commute aux changements de base quelconques et que
$\mathcal D/H \to S$ est une courbe propre lisse (voir aussi
\cite{BertinMezardQuotient} pour une autre preuve). Par suite, le
genre du quotient $\mathcal D/H$ est localement constant. Il est donc
naturel de s'intéresser aux champs algébriques $\mathcal M_{g, g'}[H
\subset G]$  des courbes propres et lisses de genre $g$ telles que le
quotient par $H$ est de genre $g'$. On obtient ainsi une décomposition  
$$\mathcal M_g[H \subset G]= \coprod_{g'} \mathcal M_{g, g'}[H \subset G]$$
en sous champs ouverts et fermés.

Dans la suite, nous fixerons un sous-groupe distingué $H \lhd G$ et
noterons $\mathcal M_{g, g'}[H \lhd G]$ afin de conserver visuellement
l'information. 

La formation du quotient étant ici fonctorielle, on obtient un
morphisme $\mathcal M_{g, g'}[H \lhd G] \to \mathcal M_{g'}[G/H]$
définit par $\mathcal D \to \mathcal D/H$.

Supposons que $\mathcal D/H$ est de genre $g'$ et notons $B$ le diviseur
de branchement du morphisme quotient. On montre que c'est un diviseur
de Cartier relatif sur $\mathcal D/H$ de degré  $n=|H|((2g-2)-(2g'-2)
|H|)$ te dont la formation commute aux changements de base. 

Notons $\mathcal M_{g'}^{[n]}[G/H]$ l'espace des modules de genre $g'$
munie d'une action de $G/H$ fidèle dans chaque fibre, et d'un diviseur
de Cartier relatif de degré $n$ (cf. \cite{Romagny} pour l'algébricité
de ce champ).

D'après ce qui a été dit précédemment on a un morphisme 
$$\begin{array}{cccc}
\Phi \colon &\mathcal M_{g, g'}[H \lhd G] & \to & \mathcal
M_{g'}^{[n]}[G/H] \\ 
&(\mathcal D \to S) & \mapsto & (\mathcal D/H \to S, B).
\end{array}$$

L'espace tangent relatif de ce morphisme possède un interprétation
simple en termes de déformations infinitésimales. Pour montrer cela,
on utilise une généralisation du principe local-global (infinitésimal)
de Bertin et Mézard obtenu dans \cite{BertinMezard}. En effet, soit
$k$ un corps, $D \to \Spec k$  une courbe propre et lisse munie d'une
action fidèle de $G$ et telle que le quotient $C:=D/H$ soit de genre
$g'$. Notons $\pi\colon D \to C$ le morphisme quotient et $\H^i_G(D,
\Omega^\vee_{D/k})$ la cohomologie équivariantes du faisceau
tangent. On a une suite exacte provenant de la suite spectrale
associée au foncteur dérivé de $\mathcal F \mapsto \H^0_{G/H}(D/H,
\pi_*^H \mathcal F)=\H^0_G(D, \mathcal F)$~: 
\begin{multline*}
0 \to \H^1_{G/H}(D/H, \pi_*^H \Omega_{D/k}^\vee) \to \H^1_G(D,
\Omega^\vee_{D/k}) \to \H^0_{G/H}(D/H, \R^1 \pi_*^H \Omega_{D/k}^\vee)
\to \\ 
\H^2_{G/H}(D/H, \pi_*^H\Omega_{D/k}^\vee) \to \H^2_G(D,
\Omega^\vee_{D/k}). 
\end{multline*}
Cette suite exacte a une interprétation très simple en termes de
déformations. Le terme $\H^1_G(D, \Omega^\vee_{D/k})$ classifie les
déformations $G$-équivariantes du premier ordre de $D$. Le terme
$\H^0_{G/H}(D/H, \R^1 \pi_*^H \Omega_{D/k}^\vee)$ classifie certaines
déformations $H$-équivariantes locales dont l'image doit être nulle
dans $\H^2_{G/H}(D/H, \pi_*^H\Omega_{D/k}^\vee)$ afin de définir une
déformation globale. 

Les termes $\H^i_{G/H}(D/H, \pi_*^H\Omega_{D/k}^\vee)$ peuvent
s'interpréter en terme de déformations $G/H$-équivariantes de $C$
munie de son diviseur $B$. En effet, on a un morphisme canonique de
faisceau $\pi^H_* \Omega_{D/k}^\vee \to \Omega_{C/k}^\vee$ qui est
injectif et dont le conoyau $\O_B$ correspond au diviseur de
branchement (avec multiplicité). 

On a donc une suite exacte longue de cohomologie
\begin{multline*}
\H^0_{G/H}(C, \O_B) \to \H^1_{G/H}(D/H, \pi^H_* \Omega_{D/k}^\vee) \to
\H^1_{G/H}(C, \Omega_{C/k}^\vee) \\ 
\to \H^1_{G/H}(C, \O_B) \to \H^2_{G/H}(D/H, \pi^H_* \Omega_{D/k}^\vee)
\to  \H^2_{G/H}(C, \Omega_{C/k}^\vee) 
\end{multline*}
dans laquelle le terme $ \H^1_{G/H}(C, \Omega_{C/k}^\vee)$ correspond
aux déformations $G/H$-équivariantes  de $C$ et $\H^2_{G/H}(C,
\Omega_{C/k}^\vee)$ correspond aux obstructions à la déformations
$G/H$-équivariantes de $C$. 

Finalement, le terme $\H^0_{G/H}(C, \O_B)$ correspond aux déformations
équivariantes du diviseur de branchement.  
%% On remarque que si $C$ a peu d'automorphismes $G/H$-équivariant
%% (e.g. g' \ge 2) alors le premier morphisme est injectif.

On peut donc préciser l'espace tangent relatif du morphisme $\Phi$.

\begin{proposition}
\label{IdentificationEspaceTangent}
L'espace tangent relatif du morphisme $\mathcal M_{g, g'}[H \lhd G]
\to \mathcal M_{g'}^{[n]}[G/H]$ en un point correspondant à une courbe
$D$ s'identifie naturellement avec le noyau du
morphisme $$\H^0_{G/H}(D/H, \R^1 \pi_*^H \Omega_{D/k}^\vee) \to
\H^2_{G/H}(D/H, \pi_*^H\Omega_{D/k}^\vee),$$  
l'espace source classifiant les déformations locales $H$-équivariantes
qui sont invariantes sous l'action de $G/H$. 
\end{proposition}

En particulier, les déformations induites de $B$ et de $C$ sont nulles
dans l'espace tangent relatif.

\subsection{Lien avec un espace de torseur et rigidification}

Considérons une courbe propre et lisse $\mathcal D \to S$ munie d'une
action de $G$ fidèle dans chaque fibre et notons $B$ le diviseur de
branchement du morphisme $\pi\colon\mathcal D \to \mathcal D/H$. Alors
par définition le morphisme $\mathcal D \setminus |\pi^{-1}(B)| \to
\mathcal D/H \setminus |B|$ est étale et est donc un $H$-torseur qui
est invariant sous l'action de $G/H$ d'après la proposition
\ref{GalGalGal}.  

Notons $h:\mathfrak C \to \mathcal M_{g'}^{[n]}[G/H]$ la courbe
universelle, $\mathfrak B$ le diviseur universel sur cette courbe et
$\mathfrak f\colon\mathfrak U= \mathfrak C \setminus |\mathfrak B|
\to \mathcal M_{g'}^{[n]}[G/H]$.  

\medskip

La construction ci-dessus définit un morphisme naturel
$$\mathcal M_{g, g'}[H \lhd G] \to [\R^1 \mathfrak f_* H]^{G/H}.$$
D'autre part, comme $H$ agit trivialement sur $[\R^1 \mathfrak f_* H]$
par construction, on a une factorisation 
$$\xymatrix{
\mathcal M_{g, g'}[H \lhd G] \ar[r] \ar[rd] & \mathcal M_{g, g'}[H
\lhd G] /\hskip -1mm/ H \ar[d] \\ 
& [\R^1 \mathfrak f_* H]^{G/H}
}$$
où $\mathcal M_{g, g'}[H\lhd G] /\hskip -1mm/ H$ désigne le
$2$-quotient obtenu en quotientant les schémas d'automorphismes par le
sous-groupe distingué $H$ (cf. \cite{Romagny} \S I.3 pour plus de
précisions).

Si on regarde la restriction de ce morphisme à la caractéristique $0$
alors on sait que $[\R^1 \mathfrak f_* H]^{G/H}$ est représentable par
un champ algébrique et on peut même montrer que le morphisme construit
ci-dessus est une immersion fermée.  Nous y reviendrons dans le
paragraphe \S \ref{CalculModere}. 

Notons $[\R^1 \mathfrak f_* H]^{G/H}_G$ le sous-champs de $[\R^1
\mathfrak f_* H]^{G/H}$ classifiant les $H$-torseurs $V\to U$
invariants sous $G/H$, qui sont géométriquement intègre sur $S$ et tel
que le groupe des automorphismes de $V \to U/(G/H)$ est $G$.

\begin{theorem}
Il existe une altération $X \to \left([\R^1 \mathfrak f_*
  H]^{G/H}_G\right)_{\gen = g, red}$ avec $X$ réduit et un diagramme
commutatif 
$$\xymatrix{ & X \ar[ld] \ar[d] \\
\mathcal M_{g, g'}[H \lhd G]_{red} \ar[r] & \left([\R^1 \mathfrak f_*
  H]^{G/H}_G\right)_{\gen = g, red}
}$$
où $\mathfrak g$ désigne le genre du torseur universel.
\end{theorem}

\preuve Oubliant l'action de $G$, c'est un corollaire de
\ref{ExtensionTorseur}. Il s'agit alors de construire l'action de
$G$. Toutefois, d'après la proposition \ref{GalGalGal}, l'action de
$G$ existe sur les torseurs. On peut alors l'étendre grace au
corollaire \ref{ProlongementMorphisme}. \qed

Il est possible d'obtenir plus d'information sur ce morphisme dans
certains cas particuliers. Pour fixer les idées, nous allons supposer
que $H=G=\ZZ/p\ZZ$, le cas $H$ plus général semblant difficile à
atteindre. 

\begin{theorem}
\label{ConstructionPrincipale}
Reprenons les notations ci-dessus et supposons que $G=H=\ZZ/p\ZZ$ avec
$p$ premier. Le morphisme canonique  
$$\mathscr{T}_{\lambda^p}(\mathfrak f) \tilde \times_{[\R^1 \mathfrak
  f_* \ZZ/p\ZZ]} \left(\mathcal M_{g, g'}[\ZZ/p\ZZ]/\hskip
  -1mm/H\right) \to \mathscr{T}_{\lambda^p}(\mathfrak f)$$ 
induit une immersion fermée dont l'image est dans 
$\left(\mathscr{T}_{\lambda^p}(\mathfrak f)\right)_{\gen=g}$.
\end{theorem}

\preuve Considérons le morphisme 
$$ \mathscr{T}_{\lambda^p}(\mathfrak f) \tilde \times_{[\R^1 \mathfrak f_*
  \ZZ/p\ZZ]} \left(\mathcal M_{g, g'}[\ZZ/p\ZZ]/\hskip -1mm/\ZZ/p\ZZ\right)
  \to \mathscr{T}_{\lambda^p}(\mathfrak f).$$ 
Par définition, un point $\Spec k \to \left(\mathcal M_{g,
    g'}[\ZZ/p\ZZ]/\hskip -1mm/\ZZ/p\ZZ\right)$ induit un
$\ZZ/p\ZZ$-torseur dont le genre est $g$, et également un
$G$-torseur. On obtient donc un morphisme  
$$ \left(\mathscr{T}_{\lambda^p}(\mathfrak f) \tilde \times_{[\R^1 \mathfrak f_*
    \ZZ/p\ZZ]} \left(\mathcal M_{g, g'}[\ZZ/p\ZZ]/\hskip
    -1mm/\ZZ/p\ZZ\right)\right)_{red} \to
\mathscr{T}_{\lambda^p}(\mathfrak f)_{\gen=g}$$ 
 
Pour montrer que ce morphisme est une immersion fermée, nous allons
vérifier les hypothèses de la proposition \ref{RelevementImpliqueIso}
(on prendra garde que les champs en présence sont en fait des limites
de champs noethériens d'après la proposition
\ref{RepresentabiliteLimite}).  

Le premier point s'obtient en prenant la normalisation des torseurs
définis sur un corps, ce qui donne automatiquement une courbe lisse de
bon genre car on ne considère que des corps algébriquement clos. 

Le point concernant les anneaux de valuations discrètes s'obtient en
fait de la même manière modulo une altération \emph{via} le corollaire
\ref{ExtensionTorseur} qui montre qu'on obtient bien une
courbe lisse (voir aussi \cite{Oort_Sekiguchi_Suwa}, Lemma IV.2.3).

Il reste donc à voir la nullité de l'espace tangent relatif. Pour
cela, on peut se placer au-dessus d'un point géométrique de $\mathcal
M_{g'}^{[n]}$ de corps résiduel $k$.  

Supposons que $k$ est de caractéristique première à $p$. La
description explicite de la proposition
\ref{IdentificationEspaceTangent} permet de voir que l'espace tangent
relatif du morphisme 
$$\mathcal M_{g, g'}[\ZZ/p\ZZ] \to \mathcal M_{g'}^{[n]}$$ 
en le point considéré est nul car le terme $\R^1 \pi_*^H$ est nul. Le
résultat est donc évident. 

Supposons maintenant que $k$ est de caractéristique $p$. Suivant
Harbater dans \cite{Harbater}, on peut identifier
$\mathscr{T}_{\lambda^p}(f)$ à $\frac{1}{x}k[\frac{1}{x}]$.  

D'autre part, on peut remplacer l'espace tangent en un point de
$\mathcal M_{g, g'}[\ZZ/p\ZZ]$ par la description explicite obtenue
dans la proposition \ref{IdentificationEspaceTangent}. 

Finalement, le résultat provient de l'étude des déformations non
obstruées obtenue par Pries dans \cite{PriesAJM} Theorem 2.2.10
(description de $\mathcal M_{\phi}$). \qed

\begin{remark}
Un analogue dans le cas $H\not =\ZZ/p\ZZ$ semble difficile à obtenir
car l'étude de l'espace tangent est alors plus compliquée (le groupe
$H$ n'étant même pas \emph{a priori} cyclique) et \cite{PriesAJM} ne
couvre que le cas $\ZZ/p\ZZ \rtimes \mmu_\ell$. 
\end{remark}

\section{Espace des modules dans le cas modéré}
\label{CalculModere}

Nous souhaitons appliquer les résultats des sections précédentes au
cas modéré afin d'obtenir une description simple des espaces $\mathcal  
M_{g}[G]$. Pour cela, il convient tout d'abord de préciser quelques
propriétés des actions modérées. La première concerne le diviseur de
branchement, qu'on peut en fait remplacer fonctoriellement par un
diviseur étale. 

\begin{lemma}
Soit $S$ un schéma, $G$ un groupe fini dont l'ordre est inversible sur
$S$, $\mathcal D \to S$ une courbe propre et lisse munie d'une action
de $G$ fidèle dans chaque fibre. Notons $\pi\colon \mathcal D \to
\mathcal D/G$ le morphisme quotient. Alors le conoyau de
l'homomorphisme $\pi_*^G \Omega_{\mathcal D/S}^\vee \to
\Omega_{\mathcal C/S}^\vee$ définit un diviseur de Cartier étale sur
$S$ et sa formation commute au changement de base.
\end{lemma}

\preuve Comme $G$ est réductif sur $S$, la formation de $\pi_*^G
\Omega_{\mathcal D/S}^\vee$ commute aux changements de base. Il
s'ensuit que le conoyau est plat. Il s'agit alors de montrer qu'il
définit un schéma étale sur $S$ mais ceci peut se faire dans le cas où
$S$ est le spectre d'un corps algébriquement clos. Le calcul est alors
direct. \qed

Nous noterons $B_r$ le diviseur ainsi obtenu, et le nommerons le
diviseur de branchement \emph{réduit}. En effet, un calcul explicite
montre qu'il a même schéma sous-jacent que le diviseur de branchement
classique, ce dernier étant en fait une puissance de $B_r$.

Donnons nous donc une courbe propre et lisse $h\colon \mathcal C \to S$
ainsi qu'une immersion ouverte $j\colon\mathcal U \to \mathcal C$ qui
est le complémentaire de $\ell> 0$ sections disjointes.

Plaçons nous dans le cas où $G= (\ZZ/p\ZZ)$ pour un nombre premier
$p$, et $S$ est un $\Spec \ZZ[\zeta, \frac{1}{p}]$-champ, $\zeta$
désignant une racine primitive $p$-ième de l'unité. Il est alors
possible de calculer complètement la plupart des termes de la suite
\eqref{PrincipeLocalGlobal}. 

En effet, par dualité de Poincaré (cf. \cite{Milne} Proposition
VI.11.8) on a un isomorphisme canonique 
$$\R^2 h_* (\ZZ/p\ZZ) = \ZZ/p\ZZ(-1).$$

Si on choisit un isomorphisme $\mmu_p \to \ZZ/p\ZZ$ on peut construire
un isomorphisme $h_* \R^1 j_* (\ZZ/p\ZZ) \to (\ZZ/p\ZZ)^\ell$ qui ne
dépend que du choix d'un ordre sur les sections et de l'isomorphisme
$\ZZ/p\ZZ \to \mmu_p$. On peut donc réécrire cette suite exacte sous
la forme 

$$0 \to \R^1 h_* \ZZ/p\ZZ \to \R^1 f_*\ZZ/p\ZZ \to (\ZZ/p\ZZ(-1))^\ell
\to \ZZ/p\ZZ(-1) \to 0$$ 
et cette suite ne dépend que de paramètres explicites.

\medskip

Nous allons maintenant appliquer les résultats obtenus jusqu'à
maintenant à l'étude de l'espace $\mathcal M_{g, g'}[\ZZ/p\ZZ] \times
\Spec \ZZ[\zeta, \frac{1}{p}]$.  On sait que ce champ algébrique est
lisse sur $ \Spec \ZZ[\zeta, \frac{1}{p}]$, il est donc réduit. 

Pour étudier cet espace, il est commode d'en introduire un autre qui
est un peu plus rigide. Notons $\ell:= \frac{2g-2-p(2g'-2)}{p-1}$. On
sait que si $\mathcal M_{g, g'}[\ZZ/p\ZZ]$ est non vide, alors d'après
la formule de Hurwitz ce nombre est un entier (la réciproque est  
également vraie). On supposera donc que c'est le cas.

Notons $\mathcal M_{g, g'}[\ZZ/p\ZZ]_\ell$ l'espace des modules des
couples $(\mathcal C \to S, (e_1, \ldots, e_\ell))$ où $\mathcal C \to
S$ est une courbe propre et lisse de genre $g$ munie d'une action de
$\ZZ/p\ZZ$ fidèle dans chaque fibre, et les $e_i$ sont des sections
disjointes $S \to \mathcal C$ invariantes sous $\ZZ/p\ZZ$ telles que le lieu
de ramification est précisément $\cup_i { \supp(e_i)}$. En oubliant
les sections $e_i$, on obtient un morphisme 
 $$\mathcal M_{g, g'}[\ZZ/p\ZZ]_\ell\times \Spec \ZZ[\zeta, \frac{1}{p}]\to \mathcal M_{g,
  g'}[\ZZ/p\ZZ]\times \Spec \ZZ[\zeta, \frac{1}{p}]$$ 
qui est un revêtement étale.

D'autre part, $\mathfrak S_\ell$ agit sur ce morphisme par permutation
des $e_i$. Il est alors aisé de voir que le morphisme induit 
$$\left(\mathcal M_{g, g'}[\ZZ/p\ZZ]_\ell \times \Spec \ZZ[\zeta,
  \frac{1}{p}] \right) / \mathfrak S_\ell \to \mathcal M_{g,
  g'}[\ZZ/p\ZZ]\times \Spec \ZZ[\zeta, \frac{1}{p}]$$ 
est un isomorphisme. 

Nous allons donc nous concentrer sur $\mathcal M_{g,
  g'}[\ZZ/p\ZZ]_\ell \times \Spec \ZZ[\zeta, \frac{1}{p}]$ dans un
premier temps. 

Notons $\mathcal M_{g'}^\ell$ l'espace des modules des courbes de
genre $g'$ munie de $\ell$ sections disjointes, $\mathfrak
h:\mathfrak C \to \mathcal M_{g'}^\ell \times \Spec \ZZ[\zeta,
\frac{1}{p}]$ la courbe universelle, $e_1 \ldots, e_\ell$ les sections
universelles et $\mathfrak U$ le complémentaire dans $\mathfrak C$ du
support des $e_i$. Notons de plus  $\mathfrak f:\mathfrak U \to
\mathcal M_{g'}^\ell \times \Spec \ZZ[\zeta, \frac{1}{p}]$ le
morphisme induit. 

On sait que $[\R^1 f_* \ZZ/p\ZZ]$ est représentable d'après la
proposition \ref{RepresentabiliteCasModere}, et en utilisant les mêmes
techniques que pour la proposition \ref{ConstructionPrincipale}, on
montre aisément qu'on a un isomorphisme 
$$\left(\mathcal M_{g, g'}[\ZZ/p\ZZ]_\ell \times \Spec \ZZ[\zeta,
  \frac{1}{p}]\right) /\hskip -1mm/ (\ZZ/p\ZZ) \to [\R^1 f_*
\ZZ/p\ZZ]_{\gen=g}.$$ 

Reprenons la suite exacte \eqref{DevissageFaisceau}. Dans le cas où $p$
est premier à toute les caractéristiques résiduelle et que la base
contient une racine primitive $p$-ième de l'unité,
les calculs obtenus en début de section peuvent alors être appliqués
ici. On en déduit une suite exacte

$$0 \to [\R^1 \mathfrak h_* \ZZ/p\ZZ] \to [\R^1 \mathfrak f_*
\ZZ/p\ZZ] \to (\ZZ/p\ZZ)^\ell_{\mathcal M_{g'} \times \Spec \ZZ[\zeta,
  \frac{1}{p}]} \to \ZZ/p\ZZ(-1)_{\mathcal M_{g'} \times \Spec
  \ZZ[\zeta, \frac{1}{p}]} \to 0.$$ 

Prenant en compte la condition $\gen=g$, on voit alors que l'espace
$$\left(\mathcal M_{g, g'}[\ZZ/p\ZZ]_\ell \times \Spec \ZZ[\zeta,
  \frac{1}{p}]\right) /\hskip -1mm/(\ZZ/p\ZZ)$$
 est muni d'une action
de $[\R^1 h_* \ZZ/p\ZZ]$ et que le quotient s'identifie naturellement avec 
$$\left\{(m_i) \in (\ZZ/p\ZZ \setminus \{0\})^\ell | \sum m_i =
  0\right\}_{\mathcal M_{g'} \times \Spec \ZZ[\zeta, \frac{1}{p}]}.$$ 

Utilisant l'isomorphisme 
$$\left(\mathcal M_{g, g'}[\ZZ/p\ZZ]_\ell \times \Spec \ZZ[\zeta,
  \frac{1}{p}]\right)/\mathfrak S_\ell\to \mathcal M_{g,
  g'}[\ZZ/p\ZZ]\times \Spec \ZZ[\zeta, \frac{1}{p}]$$ 
on obtient alors la proposition suivante.

\begin{proposition}
\label{IdentificationComposanteModeree}
Le champ algébrique $\left(\mathcal M_{g, g'}[\ZZ/p\ZZ] \times \Spec
  \ZZ[\zeta, \frac{1}{p}]\right) /\hskip -1mm/(\ZZ/p\ZZ)$  est
naturellement muni d'une action de $[\R^1 h_* \ZZ/p\ZZ]$ et le
quotient s'identifie avec 
$$\left(\left\{(m_i) \in (\ZZ/p\ZZ \setminus \{0\})^\ell | \sum m_i =
    0\right\}/\mathfrak S_\ell\right)_{\mathcal M_{g'} \times \Spec
  \ZZ[\zeta, \frac{1}{p}]}$$ 

De plus, le morphisme quotient induit une bijection entre les
composantes irréductibles des espaces de départ et d'arrivé.
\end{proposition}

\preuve Seule la dernière assertion est à démontrer, elle provient de
\cite{Cornalba} Remark 1. \qed 

Comme $\mathcal M_{g'}^\ell$ est géométriquement irréductible, on a une
description complète des composantes irréductibles de $\mathcal M_{g,
  g'}[\ZZ/p\ZZ] \times \Spec \ZZ[\zeta, \frac{1}{p}]$. 

Dans le cas où $p$ n'est pas premier, il est possible d'obtenir une
description semblable. Les différences proviennent du calcul de la
fonction $\gen$ ainsi que de l'irréductibilité des composantes. 
%% l'argument de Cornalba fonctionne bien mais donne probablement un autre résultat !!!!

\end{document}